\documentclass[conference]{IEEEtran}
\IEEEoverridecommandlockouts
\usepackage{graphicx}
\usepackage{array}
\usepackage{textcomp}
\usepackage{cite}
\usepackage{url}
\usepackage{amsmath}
\usepackage{multirow}
\usepackage{bm}
\usepackage{verbatim}
\usepackage{array}
\newcolumntype{L}[1]{>{\raggedright\let\newline\\\arraybackslash\hspace{0pt}}m{#1}}
\newcolumntype{C}[1]{>{\centering\let\newline\\\arraybackslash\hspace{0pt}}m{#1}}
\newcolumntype{R}[1]{>{\raggedleft\let\newline\\\arraybackslash\hspace{0pt}}m{#1}}

\usepackage{etoolbox}
\usepackage{xcolor}
\newcommand{\splitatcommas}[1]{%
	\begingroup
	\begingroup\lccode`~=`, \lowercase{\endgroup
		\edef~{\mathchar\the\mathcode`, \penalty0 \noexpand\hspace{0pt plus 0em}}%
	}\mathcode`,="8000 #1%
	\endgroup
}

\ifCLASSOPTIONcompsoc
\usepackage[caption=false,font=normalsize,labelfont=sf,textfont=sf]{subfig}
\else
\usepackage[caption=false,font=footnotesize]{subfig}
\fi

\ifCLASSINFOpdf
\else
\fi
\begin{document}
%
\title{Community Microgrid Planning Considering Building Thermal Dynamics}

\author{\IEEEauthorblockN{Xiaohu Zhang\IEEEauthorrefmark{1},
		Desong Bian\IEEEauthorrefmark{1},
		Di Shi\IEEEauthorrefmark{1},
		Zhiwei Wang\IEEEauthorrefmark{1} and
		Guodong Liu\IEEEauthorrefmark{2}}
	\IEEEauthorblockA{\IEEEauthorrefmark{1}GEIRI North America, San Jose, CA, USA}
	\IEEEauthorblockA{\IEEEauthorrefmark{4}Oak Ridge National Laboratory, Oak Ridge, TN, USA\\
		Email:xiaohu.zhang@geirina.net}
	\thanks{This work is supported by SGCC Science and Technology Program under project Hybrid Energy Storage Management Platform for Integrated Energy System.}}



\maketitle

\begin{abstract}
In this paper, the building thermal dynamic characteristics are introduced in the community microgrid (MG) planning model. The proposed planning model is formulated as a mixed integer linear programming (MILP) which seeks to determine the optimal deployment strategy for various distributed energy resources (DER). The objective is to minimize the annualized cost of community MG, including investment cost in DER, operation cost for dispatchable fuel-generators (DFG), energy storage system (ESS) degradation cost, energy purchasing and peak demand charge at PCC, customer discomfort cost due to the room temperature deviation and load curtailment cost. Given the slow thermal dynamic characteristics of buildings, the heating, ventilation and air-conditioning (HVAC) system in the proposed model is treated as a demand side management (DSM) component, whose dispatch commands are provided by the central MG controller. Numerical results based on a community MG comprising of 20 residential buildings demonstrate the effectiveness of the proposed planning model and the benefits of introducing the building thermal dynamic model.

\end{abstract}


%
\IEEEpeerreviewmaketitle

\section*{Nomenclature}
\subsection*{Indices and Sets}
\addcontentsline{toc}{section}{Nomenclature}
\begin{IEEEdescription}[\IEEEusemathlabelsep\IEEEsetlabelwidth{$V_1,V_2,V_3, V_4$}]
	\item[$n$] Index of DFG.
	\item[$l$] Index of energy blocks offered by DFG.
	\item[$w, \ v$] Index of WT and PV.
	\item[$b$] Index of ESS.
	\item[$m,d,t$] Index of month, day, time interval.
	\item[$\Omega_{n}, \Omega_{l}$] Set of DFGs and their energy blocks.
	\item[$\Omega_{w},  \Omega_{v}, \ \Omega_{b}$] Set of WTs, PVs and ESSs.
	\item[$\Omega_m, \Omega_{d},  \Omega_{t}$] Set of months, days and time intervals.
	\item[$\Omega_{dt}^m$] Set of days and time intervals belong to month $m$.
\end{IEEEdescription}

\subsection*{Variables}
\addcontentsline{toc}{section}{Nomenclature}
\begin{IEEEdescription}[\IEEEusemathlabelsep\IEEEsetlabelwidth{$V_1,V_2,V_3$}]
	\item[$\delta_w, \ \delta_v$] Binary variable associated with installing an WT/PV.
    \item[$\delta_n, \ \delta_b$] Binary variable associated with installing a DFG/ESS.
	\item[$\alpha_{ndt}$] Binary variable indicating on/off status of DFG $n$ on day $d$ at time $t$.
	\item[$\beta_{ndt}$] Binary variable associated with the start-up status of DFG $n$ on day $d$ at time $t$.
	\item[$P_{ndtl}$] Power generation from the $l$-th block of energy offered by DFG $n$ on day $d$ at time $t$.
    \item[$P_{ndt}$] Power output of DFG $n$ on day $d$ at time $t$.
    \item[$P_{bdt}^c,\ P_{bdt}^d$] Charging/discharging power of ESS $b$ on day $d$ at time $t$.
    \item[$SOC_{bdt}$] State of charge of ESS $b$ on day $d$ at time $t$.
    \item[$P_{wdt}, \  P_{vdt}$] Power output from WT/PV $w/v$ on day $d$ at time $t$.
    \item[$\Delta P_{hdt}^{ls}$] Load shedding in house $h$ on day $d$ at time $t$.
    \item[$P_{dt}^{PCC}$] Power exchange at PCC on day $d$ at time $t$.
    \item[$P_m^{pk}$] Peak demand at PCC in month $m$.
\end{IEEEdescription} 

\subsection*{Parameters}
\addcontentsline{toc}{section}{Nomenclature}
\begin{IEEEdescription}[\IEEEusemathlabelsep\IEEEsetlabelwidth{$V_1,V_2,V_3,V_4$}]
	\item[$A_w,A_v$] Annuity factor of WT/PV.
	\item[$A_n,A_b$] Annuity factor of DFG/ESS.
   \item[$C_w^P,C_v^P,C_n^P$] Cost of power capacity of WT/PV/DFG.
   \item[$C_b^P,C_b^e$] Cost of power/energy capacity of ESS.
   \item[$a_{nl}$] Marginal cost of the $l$-th block of energy offered by DFG $n$.
   \item[$k_n,e_n$] No-load and startup cost for DFG $n$.
   \item[$\lambda_{dt}^{PCC}$] Electricity price at PCC on day $d$ at time $t$.
   \item[$c_{b}$] Degradation cost of ESS $b$.
   \item[$w_{h}$] Customer discomfort cost of house $h$.
   \item[$d^{ls}$] Penalty for load shedding.
   \item[$\lambda_m^{dc}$] Peak demand in month $m$. 
   \item[$P_n^{\min},P_n^{\max}$] Minimum/maximum output of DFG $n$.
   \item[$P_{nl}^{\max}$] The maximum output power of the $l$-th energy block for DFG $n$.
   \item[$P_b^{\max},E_b^{\max}$] Power/energy capacity of ESS $b$.
   \item[$\gamma_b^{\min},\gamma_b^{\max}$] Minimum/maximum energy level of ESS $b$.
   \item[$\eta_b^c,\eta_b^d$] Charging/discharging efficiency of ESS $b$.
   \item[$P_w^{\max},P_{v}^{\max}$] Power capacity of WT/PV.
   \item[$r_{wdt},r_{vdt}$] Available capacity factor of WT/PV $w/v$ on day $d$ at time $t$.
   \item[$T_{hdt}^d$] The desired indoor temperature of house $h$ on day $d$ at time $t$.
   \item[$\theta_{hdt}$] Allowed temperature deviation of house $h$ on day $d$ at time $t$. 
   \item[$P_{hdt}^{ls,\max}$] Maximum load curtailment of house $h$ on day $d$ at time $t$. 
\end{IEEEdescription} 
Other symbols are defined as required in the text.

\section{Introduction}
\IEEEPARstart{A}{}microgrid (MG) is a small-scale, low-voltage active distribution network consisting of different distributed generators (DG), energy storage system (ESS) and responsive loads, which is operated either in grid-connected or islanded mode \cite{mybibb:nikos_MG,mybibb:certs_concept}. In grid-connected mode, a MG can be regarded as a controllable entity as it can not only export/import power from the main grid, but also offer various ancillary services such as voltage and frequency regulation through the point of common coupling (PCC). In islanded mode, a MG has the capability of satisfying loads locally, which enhances energy reliability and resiliency.  Thus,  there is a great interest in the implementation and deployment of MG from both industry and academia \cite{mybibb:yongli,mybibb:jiajun,mybibb:zhehan1}.

The objective of MG expansion planning is to determine the optimal size and type of distributed energy resources (DER), i.e., DG and ESS. A rationally planned MG should ensure the economic benefits of MG deployments and further justify the return on investments (ROI). The models and solution approaches for the MG planning have been extensively studied in the technical literature. In \cite{mybibb:ess_allocation_GA}, the optimal allocation and economic operation of ESS in MG are addressed by genetic algorithm (GA). The authors in \cite{mybibb:chen_pso_MG} leverage particle swarm optimization (PSO) to determine the optimal size of distributed energy resources (DER) in a community MG. The obtained results are compared with HOMER software \cite{mybibb:homer_web}. To improve the voltage profile and reduce power losses, reference  \cite{mybibb:dg_sensitivity_1} proposes an approach to optimally place DG in distribution network based on continuous power flow (CPF).

With the advances in branch-and-bound algorithm, mathematical programming has also been widely employed to solve the planning problem. Reference \cite{mybibb:yishen} proposes a DER sizing model in a hybrid AC/DC MG configuration. The authors in \cite{mybibb:reliability_constrained_mg_planning} propose a method to calculate the optimal ESS size in a MG considering reliability criterion. The complete model is formulated as a mixed integer linear programming (MILP) problem and solved by commercial solvers. The co-planning of renewable energy resources (RES) and ESS in a residential MG is presented in \cite{mybibb:size_re_ess}. The intermittent characteristics of renewable generation are captured by stochastic programming. In \cite{mybibb:microgrid_planning_uncertain}, the MG planning model considering different uncertainties such as forecast errors for loads, variable renewable generation and islanding incidents is proposed. To relieve the computational burden, the complete optimization model is decomposed into the investment master problem and operation subproblems. The authors in \cite{mybibb:zhaoyu_dg} propose a two-stage robust optimization model to identify the optimal location of DG in MG. The column-and-constraint generation (CCG) is leveraged to solve the problem. A chance constrained information-gap decision model for multi-stage MG planning model is presented in \cite{mybibb:chance_microgrid}. A customized bilinear Benders decomposition algorithm is developed to solve the model. 

Few of the previously mentioned work have incorporated the heating, ventilation, and air conditioning (HVAC) systems and building thermal dynamics into the MG planning model. According to \cite{mybibb:vietnam_2}, the HVAC system can be treated as a promising candidate for the demand side management (DSM) given the slow thermal dynamic characteristics of buildings. Specifically, with a user preferred indoor temperature, the HVAC system can precool/preheat buildings when the electricity price is low and be switched off when the price is high without sacrificing the user comfort level. Therefore, a more useful MG investment plan can be achieved if the building thermal dynamic model is considered. \textcolor{black}{Reference \cite{mybibb:gd_size_ess} integrates the thermal dynamic model of buildings in the ESS sizing problem. To the best of the authors’ knowledge, none of the existing literatures evaluates the impacts of building thermal dynamics on MG planning model.}

This paper proposes a community MG planning model incorporating building thermal dynamics. Consider a single target year, the optimization model jointly minimizes the investment in different types of DERs, as well as the expected operation cost. The complete model is originally nonlinear and recasted as an MILP model by a linearization technique. The contribution of this paper are twofold:
\begin{itemize}
	\item \textcolor{black}{to develop a MG planning model integrating the building thermal dynamics;}
	\item \textcolor{black}{to evaluate the benefits brought by the building thermal dynamics in the MG planning process based on detailed numerical results.}
\end{itemize}

The remaining sections are organized as follows. In Section \ref{system_model}, the model of community MG and building thermal dynamics are described. Section \ref{formulation_1} illustrates details of the MG planning model. In Section \ref{results}, case studies and numerical results are provided. Finally, conclusions and future work are given in Section \ref{conclusion_future}. 

\section{System Model}
\label{system_model}
\subsection{Community MG}
\label{community_mg}
In this work, the candidate DER for the considered community MG includes dispatchable fuel-generators (DFG), RES and ESS. There also exists a number of residential houses. We assume that the load for each house is divided into HVAC load and non-HVAC load. Moreover, each house is equipped with a house energy management system (HEMS), which collects the house information, controls the house appliances and communicates with the central MG controller. With the information provided by each HEMS, forecasted renewable generations and electricity market price, the community MG controller determines the optimal scheduling of DG and ESS, the amount of power exchange at PCC and the operation states of HVAC systems.  

\subsection{HVAC System and Building Thermal Dynamics}
Traditionally, the HVAC system in a house is controlled by a temperature sensor and a controlled relay circuit. With the allowable indoor temperature settings provided by customers, the temperature sensor detects the current indoor temperature and the controlled relay circuit decides the on/off status of HVAC. As an example, in the heating mode, the controlled relay will turn the HVAC on if the indoor temperature is lower than the floor of the allowed temperature. The HVAC will be switched off until the indoor temperature reaches the ceiling of the allowed temperature. The traditional HVAC control method is referred as `simple' control in this work.

As mentioned in Section \ref{community_mg}, we assume that the central MG controller implements surrogate control of HVAC systems and offers the dispatch commands based on the scheduling model. To achieve that, the building thermal dynamics should be included in the MG scheduling framework. We leverage the third order state space model proposed in \cite{mybibb:thermal_master} to describe the building thermal dynamic characteristics as:
\begin{equation}
\bm{T}_{h,d,t+1}=\bm{A}_h\bm{T}_{hdt}+\bm{B}_h\bm{U}_{hdt}, \  \forall d, \forall h, \forall t  \label{thermal}
\end{equation}
where $\bm{T}_{dht}=[T_{hdt}^{in},T_{hdt}^{m},T_{hdt}^e]^T$ denotes the state vector and  $\bm{U}_{hdt}=[T_{dt}^{a},\Phi_{dt},(u_{hdt}^H-u_{hdt}^C)\xi_hP_h^H]^T$ represents the input vector. Specifically, $T_{hdt}^{in}$, $T_{dht}^{m}$ and $T_{dht}^{e}$ indicate the indoor temperature, the temperature of thermal accumulating layer of the inner walls and the temperature of house envelope respectively, for house $h$ on day $d$ at time interval $t$. $T_{dt}^a$ and $\Phi_{dt}$ denote the ambient temperature and solar irradiance on day $d$ at time interval $t$. Two binary variables $u_{hdt}^H$ and $u_{hdt}^C$ are introduced to flag the HVAC operation mode, i.e., heating and cooling. $\xi_h$ is the coefficient of performance (COP) of HVAC for house $h$ and $P_h^H$ denotes the rated power of HVAC system in house $h$. Finally, the coefficients matrices $\bm{A}_h$ and $\bm{B}_h$ depend on the thermal capacitance and resistance of the house, the effective window area and the fraction of solar irradiatation entering the inner walls and floor.  Reference \cite{mybibb:thermal_master,mybibb:vietnam_2} provide details regarding the derivation of building thermal dynamics, i.e., equation (\ref{thermal}).

\section{Problem Formulation}
\label{formulation_1}
In this section, the complete MG planning model is first illustrated, and then a linearization technique is leveraged to transform the nonlinear model into an MILP model.
\subsection{Optimization Model}
The complete optimization model is given by (\ref{objective})-(\ref{demand_lim}):
\begin{align}
&\min_{\Xi_{\text{OM}}} C_I+C_O \label{objective}  \\
&\text{s.t.}\ \ (\ref{thermal}) \ \  \text{and}     \nonumber \\
&C_{I}=A^w \sum_{w \in \Omega_w}C^p_wP_w^{\max}\delta_w+A^v\sum_{v \in \Omega_v}C^p_vP_v^{\max}\delta_v   \nonumber \\
&+A^n\sum_{n \in \Omega_n}C^p_nP_n^{\max}\delta_n+A^b\sum_{b\in \Omega_b}(C_b^pP_b^{\max}+C_b^eE_b^{\max})\delta_b \label{inv_exp} \\ 
&C_{O}=\pi_1\{\sum_{n \in \Omega_n}\sum_{d \in \Omega_d} \sum_{t \in \Omega_t} [\sum_{l \in \Omega_l}a_{nl}P_{ndtl}+k_n\alpha_{ndt}+e_n\beta_{ndt}]  \nonumber \\
&+\sum_{d \in \Omega_d} \sum_{t \in \Omega_t}\lambda_{dt}^{PCC}P_{dt}^{PCC}+\sum_{b \in \Omega_b}\sum_{d \in \Omega_d} \sum_{t \in \Omega_t}c_{b}(P_{bdt}^c+P_{bdt}^d)   \nonumber \\
&+\sum_{h \in \Omega_h}\sum_{d \in \Omega_d} \sum_{t \in \Omega_t}w_{h}|T_{hdt}^{in}-T_{hdt}^d|+ d^{ls}\Delta P_{hdt}^{ls}\}  \nonumber  \\
&+\pi_2\sum_{m \in \Omega_m}\lambda_{m}^{dc}P_m^{pk}  \label{oper_exp}  
\end{align}
\begin{align}
&C_I\le C_I^{\max} \label{budget}   \\
&P_{ndt}=\sum_{l \in \Omega_l}P_{ndtl}+\alpha_{ndt}P_{n}^{\min},\ \ \forall n, \forall t, \forall d \label{dg_output}  \\
&0\le P_{ndtl} \le P_{nl}^{\max}, \ \ \forall n, \forall d, \forall t, \forall l \label{dg_segment}  \\
&0\le P_{ndt} \le P_{n}^{\max}\alpha_{ndt}, \ \ \forall n, \forall d, \forall t \label{dg_output_lim}  \\
&\beta_{ndt} \ge \alpha_{ndt}-\alpha_{nd,t-1}, \ \ \forall n, \forall d, \forall t \label{st_logic}  \\
&\alpha_{ndt}\le \delta_n, \ \ \forall n, \forall d, \forall t \label{install_logic} \\
&0\le P_{bdt}^c \le \delta_bP_b^{\max}, \ \ \forall b, \forall d, \forall t \label{ch_lim} \\
&0\le P_{bdt}^d \le \delta_bP_b^{\max}, \ \ \forall b, \forall d, \forall t \label{dch_lim} \\ 
&\delta_b\gamma_{b}^{\min}E_b^{\max} \le SOC_{bdt} \le \delta_b\gamma_{b}^{\max}E_b^{\max}, \ \ \forall b, \forall d, \forall t \label{soc_lim} \\
&SOC_{bdt}=SOC_{bd,t-1}+P_{bdt}^c\eta_b^c\Delta t-P_{bdt}^d/\eta_b^d\Delta t,     \nonumber \\
&\ \ \ \ \ \ \ \ \ \  \ \ \ \ \ \ \ \ \ \ \ \ \ \ \ \ \ \ \ \ \ \ \ \ \ \ \ \ \ \ \ \ \ \ \ \ \ \ \ \ \  \forall b, \forall d, \forall t  \label{soc_cons} \\
&0\le P_{wdt} \le \delta_w P_{w}^{\max}r_{wdt}, \ \ \forall w, \forall d, \forall t  \label{wind_lim}\\
&0 \le P_{vdt} \le \delta_v P_v^{\max}r_{vdt} , \ \ \forall v, \forall d, \forall t \label{pv_lim}\\
&T_{hdt}^d-\theta_{hdt} \le T_{hdt}^{in} \le T_{hdt}^d+\theta_{hdt}, \ \ \forall h, \forall d, \forall t \label{indoor_lim}  \\
&0 \le \Delta P_{hdt}^{ls}\le P_{hdt}^{ls,\max}, \ \ \forall h, \forall d, \forall t \label{LS_lim}  \\
&\sum_{w \in \Omega_w}P_{wdt}+\sum_{v \in \Omega_v}P_{vdt}+\sum_{n \in \Omega_n}P_{ndt}+\sum_{b \in \Omega_b}(P_{bdt}^d-P_{bdt}^c) \nonumber \\
&+P_{dt}^{PCC}=\sum_{h \in \Omega_h}[(u_{hdt}^H+u_{hdt}^C)P_h^H+(P_{hdt}^O-\Delta P_{hdt}^{ls})],   \nonumber \\
&\ \ \ \ \ \ \ \ \ \ \ \ \ \ \ \ \ \ \ \ \ \ \ \ \ \ \ \ \ \ \ \ \ \ \ \ \ \ \ \ \ \ \ \ \ \ \ \ \ \ \ \ \ \ \  \forall d, \forall t  \label{power_ba} \\
&-P_{dt}^{PCC,\max} \le P_{dt}^{PCC} \le P_{dt}^{PCC,\max}, \forall d, \forall t \label{pcc_lim}  \\
&P_m^{pk} \ge P_{dt}^{PCC}, \ \ \forall m, \forall d \in \Omega_{dt}^m, \forall t \in \Omega_{dt}^m \label{demand_lim}  
\end{align}

The optimization variables of the proposed MG planning model are those in set $\splitatcommas{\Xi_{\text{OM}}=\{\delta_w,\delta_v,\delta_n,\delta_b,P_{ndtl},P_{ndt},\alpha_{ndt},\beta_{ndt},P_{bdt}^c,P_{bdt}^d,SOC_{bdt},P_{wdt},P_{vdt},u_{hdt}^H,u_{hdt}^C,\Delta P_{hdt}^{ls},T_{hdt}^{in},T_{hdt}^m,T_{hdt}^e,P_{dt}^{PCC},P_m^{pk}\}}$. The objective function seeks to minimize the total annualized cost including both the investment ($C_I$) and the operation cost ($C_O$), whose expressions are described by constraint (\ref{inv_exp}) and (\ref{oper_exp}). Specifically, the first line in (\ref{inv_exp}) denotes the investment in wind turbines and PV panels and the second line represents the investment in DFG and ESS. With respect to the annualized operation cost, i.e., (\ref{oper_exp}), several representing days are selected. The first line in (\ref{oper_exp}) indicates the operation cost for DFG which includes three terms. The first term denotes the fuel cost which is  described by the piecewise linear function. The no-load and startup cost are represented by the second and third terms. The second line denotes the energy purchasing/selling cost/benefits at the PCC of distribution network and the degradation cost of ESS. The customer discomfort cost due to the difference between indoor and desired temperature, and the involuntary load shedding cost are described by the third line. Note that the aforementioned operation costs are on daily basis and should be scaled to annualized cost by the scaling factor $\pi_1$. Finally, the last line in (\ref{oper_exp}) represents the monthly peak demand charge, which is scaled to the annualized cost by a different scaling factor $\pi_2$.

The budget on the annualized investment in DER is imposed in constraint (\ref{budget}). Constraint (\ref{dg_output}) and (\ref{dg_segment}) represent the output for DFG and ensure that they deliver at least minimum power when committed. Constraint (\ref{dg_output_lim}) enforces the generation of dispachable units to be zero if not committed. The relationship between start up indicator and on/off status of dispatchable units is illustrated in constraint (\ref{st_logic}) \cite{mybibb:scuc_classic}. Constraint (\ref{install_logic}) enforces that the DFG will produce no power when it is not installed. Constraint (\ref{ch_lim}) and (\ref{dch_lim}) denote the limits of charging and discharging power of ESS. The limit on the state of charge (SOC) is enforced by (\ref{soc_lim}). The SOC transition is illustrated by constraint (\ref{soc_cons}). Constraint (\ref{wind_lim}) and (\ref{pv_lim}) denote that the output of RES should be bounded by its available capacity. The indoor temperature is constrained by (\ref{indoor_lim}). Constraint (\ref{LS_lim}) places an upper limit on the load curtailment amount. The power balance equation is enforced by constraint (\ref{power_ba}). The limits on the amount of power exchanging with the power grid are denoted by constraint (\ref{pcc_lim}). Finally, the peak power at PCC is represented by constraint (\ref{demand_lim}). 
\subsection{Linearization}
As observed from the optimization model, the only nonlinear term is the absolute term in (\ref{oper_exp}). Fortunately, this term exists in the objective function. We introduce two positive slack variables and one additional constraint to linearize the term as follow \cite{mybibb:zheming,mybibb:svc_gm_2017}:
\begin{align}
&T_{hdt}^{in}-T_{hdt}^d+s_{hdt,1}-s_{hdt,2}=0   \label{abs_cons1}  \\
&s_{hdt,1} \ge 0, s_{hdt,2} \ge 0   \label{slack_con1}
\end{align}  
Then the absolute values term becomes:
\begin{equation}
w_{h}|T_{hdt}^{in}-T_{hdt}^d| \rightarrow w_{hdt}(s_{hdt,1}+s_{hdt,2})
\end{equation}

The complete model is transformed into an MILP, which can be efficiently solved by off-the-shelf solvers.

\section{Numerical Results}
\label{results}
We test our proposed planning model on a community MG consisting of 20 residential buildings with corresponding HVAC and non-HVAC loads. Following \cite{mybibb:gd_size_ess}, we select three days, i.e., Jan. 1st, Apr. 15th and Aug. 1st in 2016, which represent a typical winter, spring and summer day in the southeast states of US. The time resolution is chosen to be 15 minutes. More days can be easily included but the trade-off between computational efficiency and accuracy should be considered. The solar irradiance and temperature on the selected days are extracted from \cite{mybibb:data_ornl}. The peak value of non-HVAC load is 213.47 kW. The non-HVAC load profile, energy price at PCC are the same as \cite{mybibb:guodong_data}. The peak demand charge is selected to be \$2/kW for summer and \$1/kW for winter and spring. The load shedding cost is selected to be \$10/kWh and the upper bound of the load curtailment is 10\% of the non-HVAC load. The rated power ($P_h^H$) and COP ($\xi_h$) of HVAC system are selected to be 5 kW and 4 respectively. The desired indoor temperature is selected to be 21$^\circ C$ with the allowed deviation to be $\pm$2$^\circ C$. The customer discomfort cost is set to be \$0.05/$^\circ C$. Reference \cite{mybibb:thermal_master} provides the other parameters for the residential buildings. The standard error of estimate is used to represent the variety of different houses.

The candidate DER comprises wind turbines (WT), PV panels (PV), diesel engines (DE), microturbines (MT) and two types of batteries (ES1, ES2), whose investment parameters are provided in Table \ref{inv_para} and \ref{inv_ess}. The other operation parameters of DE and MT, i.e., startup cost and no-load cost, are taken from \cite{mybibb:guodong_data}. The operation cost of WT and PV are assumed to be zero. The available WT and PV capacity factors are based on \cite{mybibb:data_nrel}. Note that the annuity factors in (\ref{inv_exp}) are calculated by the interest rate and life time of DER \cite{mybibb:facts_ccg}. In this work, the interest rate is set to be 5\% and the life time of DER is 10 years.

The complete problem is implemented in YALMIP \cite{mybibb:YALMIP} and solved by CPLEX \cite{mybibb:CPLEX1}. The ``mipgap" is set to be 0.5\%. All the simulations are conducted on a computer with an Inter Core(TM) i7-6600U CPU @ 2.60 GHz and 16.0 GB RAM.

\begin{table}[!htb]
	\centering
	\caption{Investment Parameters of Candidate DG}
	\label{inv_para}
	\begin{tabular}{c c c c c c c c}
	\hline
	&Rated&Min&Capital&Num&\multicolumn{3}{c}{Fuel}  \\
	Type&Power&Output&Cost&Limit&\multicolumn{3}{c}{Cost (\$/kWh)}  \\
	&(kW)&(kW)&(\$/kW)&&$l=1$&$l=2$&$l=3$  \\
	\hline
	WT&120&0&2700&2&-&-&-  \\
	\hline
	PV&80&0&2100&2&-&-&-  \\
	\hline
	DE&60&10&540&2&0.2822&0.3732&0.4643 \\
	\hline
	MT&80&10&810&2&0.2392&0.3163&0.3936  \\
	\hline
		\end{tabular}
\end{table}   
\begin{table}[!htb]
	\centering
	\caption{Investment Parameters of Candidate ESS}
	\label{inv_ess}
	\begin{tabular}{c c c c c c c}
	\hline
	&Rated&Rated&Power&Energy&Effi-&Num   \\
	Type&Power&Energy&Cost&Cost&ciency&Limit  \\
	&(kW)&(kWh)&(\$/kW)&(\$/kWh)&&  \\
	\hline
	ES1&90&150&324&180&0.95&2 \\
	\hline
	ES2&100&200&240&216&0.85&2  \\
	\hline
	\end{tabular}
\end{table}     
We consider four cases in the simulations: 1) The budget $C_I^{\max}$ is zero for the proposed MG planning model; 2) The budget $C_I^{\max}$ is \$80k for the proposed MG planning model; 3) The budget $C_I^{\max}$ is \$100k for the proposed MG planning model; 4) The budget is the same as 3) but the simple control for HVAC system is used. 

\begin{table}[!htb]
	\centering
	\caption{MG Planning Results for Different Cases}
	\label{case_compare}
	\begin{tabular}{c c c c c}
		\hline
		&\textbf{C1}&\textbf{C2}&\textbf{C3}&\textbf{C4}  \\
		\hline
		\# of WT&0&0&1&1    \\
		\hline
		\# of PV&0&2&1&1  \\
		\hline
		\# of DE&0&1&0&0  \\
		\hline
		\# of MT&0&2&2&2  \\
		\hline
		\# of ES1&0&1&2&2  \\
		\hline
		\# of ES2&0&0&0&0  \\
		\hline
		Annualized&\multirow{2}{*}{0}&\multirow{2}{*}{77.98}&\multirow{2}{*}{98.15}&\multirow{2}{*}{98.15}  \\
		Investment Cost (\$ $10^3$)&&&&  \\
		\hline
		Annualized&\multirow{2}{*}{565.59}&\multirow{2}{*}{409.90}&\multirow{2}{*}{372.26}&\multirow{2}{*}{394.98}  \\
		Operation Cost (\$ $10^3$)&&&&  \\
		\hline
		Annualized&\multirow{2}{*}{565.59}&\multirow{2}{*}{487.88}&\multirow{2}{*}{470.42}&\multirow{2}{*}{493.14}  \\
		Total Cost (\$ $10^3$)&&&&  \\
		\hline
		Computational&\multirow{2}{*}{233.33}&\multirow{2}{*}{563.64}&\multirow{2}{*}{635.36}&\multirow{2}{*}{3.69}  \\
		Time (s)&&&&  \\
		\hline
	\end{tabular}
\end{table}

Table \ref{case_compare} provides the planning results for different cases. When the budget is zero, all the loads in the community will be covered by the purchased power from PCC. The annualized total cost is \$565.59k. For Case 2, the planning results suggest to install 2 PV, 1 DE, 2 MT and 1 ES1. The installation of DER reduces the total cost to \$487.88k despite their investment. In Case 3, 1 PV, 1 WT, 2 MT and 2 ES1 are selected and the total cost is further decreased to \$470.42k. Case 4 and Case 3 have the same investment results. Nevertheless, the operation cost for Case 3 is around  \$22.72k less than that of Case 4, which demonstrates the benefits of surrogate control of HVAC systems in our proposed model. With respect to the computation issues, Case 4 has the fastest computational time. The reason lies in that the HVAC on/off states are pre-determined by solving equation (\ref{thermal}) in the simple control so the number of binary variables in the optimization model is significantly reduced.

\begin{figure}
	\centering
	\subfloat[Simple control]{%
		\includegraphics[width=0.42\textwidth]{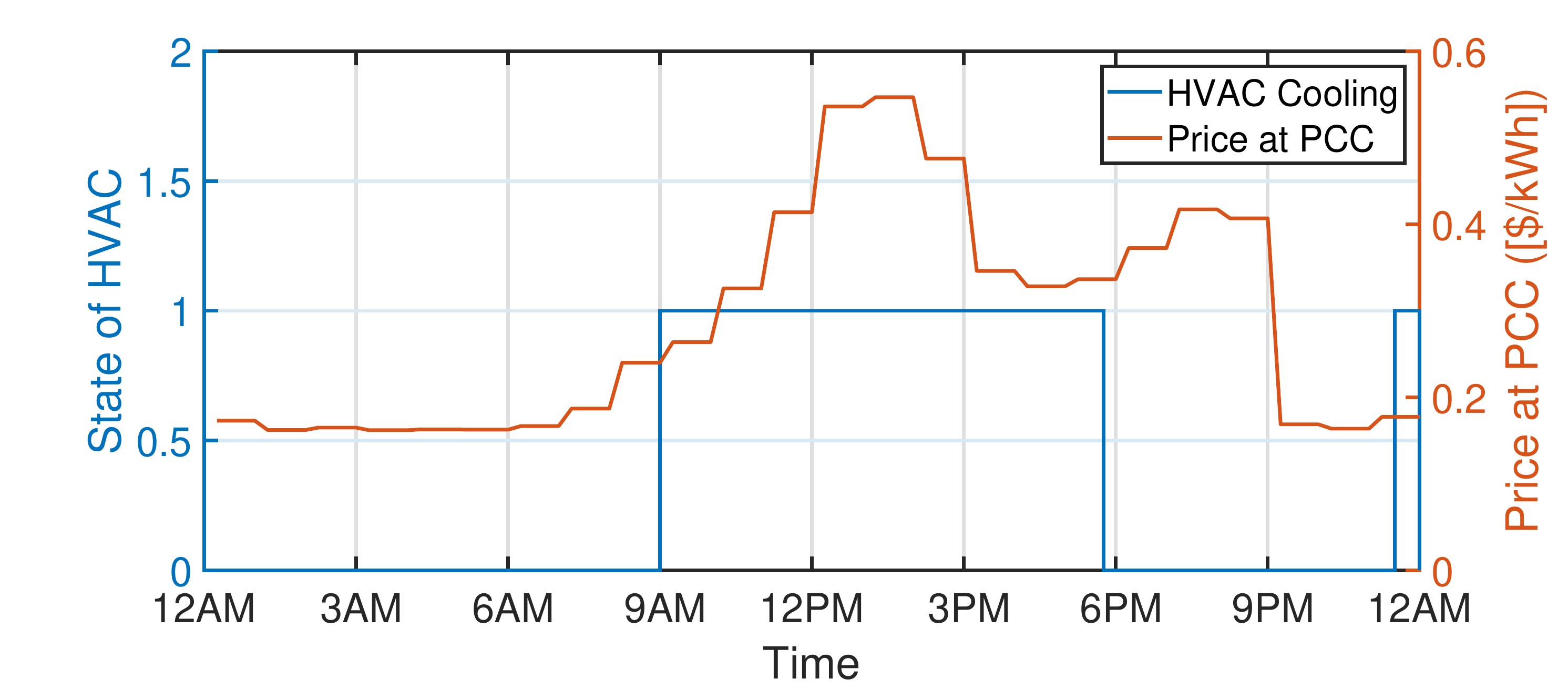}}
	\label{bus_5}
	\subfloat[Surrogate control]{%
		\includegraphics[width=0.42\textwidth]{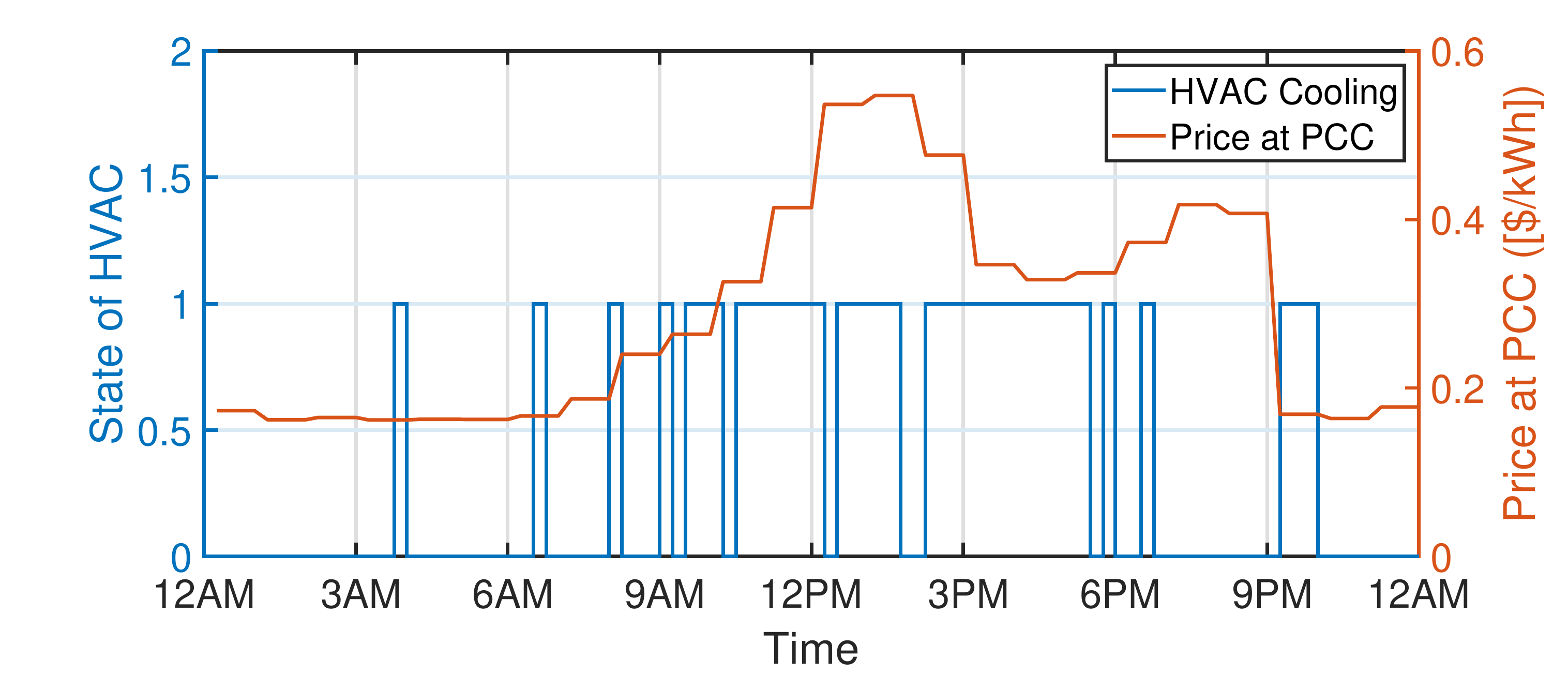}}
	\label{bus_26}
	\caption{Comparison of different control strategy of HVAC for house \#1.}
	\label{hvac_control}
\end{figure}

Fig. \ref{hvac_control} illustrates the HVAC on/off status versus the electricity price at PCC on the Summer day. As can be observed from the figure, the simple control of HVAC  does not respond to the electricity price. The HVAC is turned on to cool the house when the outdoor temperature is high, i.e, 9AM-6PM. On the other hand, the surrogate control switches on the HVAC to precool the house during the low price hours (before 9AM), and turns off the HVAC during some peak price intervals (12PM-3PM). Thus, the total operation cost is reduced.

\begin{figure}[!htb]
	\centering
	\includegraphics[width=0.45\textwidth]{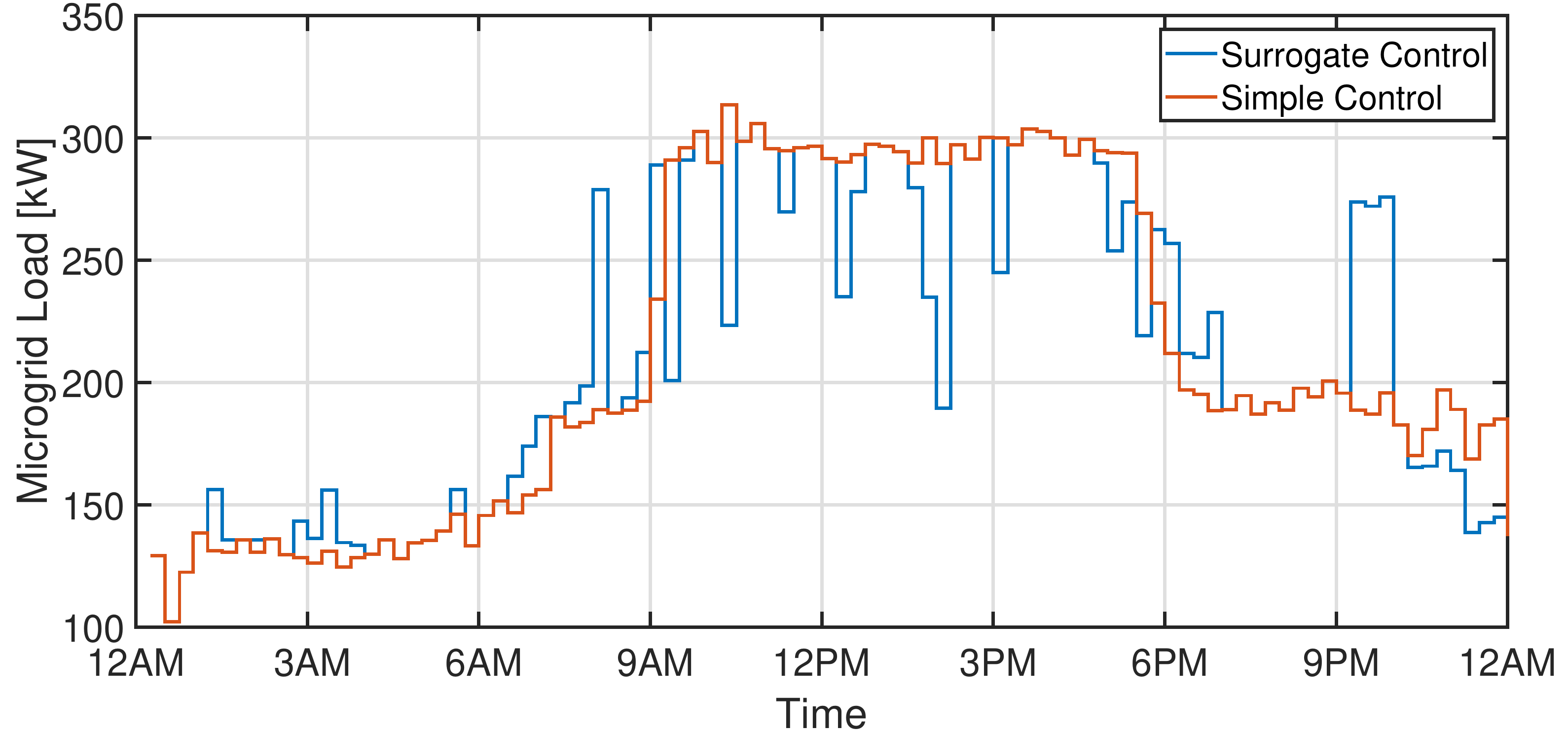}
	\caption{MG load with different control strategy on the Summer day.}
	\label{mg_load}
\end{figure}

The total load of MG by using different control strategies is depicted in Fig. \ref{mg_load}. As expected, a portion of  load is shifted to the off-peak hours, i.e., around 8AM, which is due to the precooling of the house. In addition, compared to the simple control, several off states of house HVAC lead to the MG load reduction during the peak hours (11AM-3PM), by using the surrogate control.  

\begin{figure}
	\centering
	\subfloat[Simple control]{%
		\includegraphics[width=0.45\textwidth]{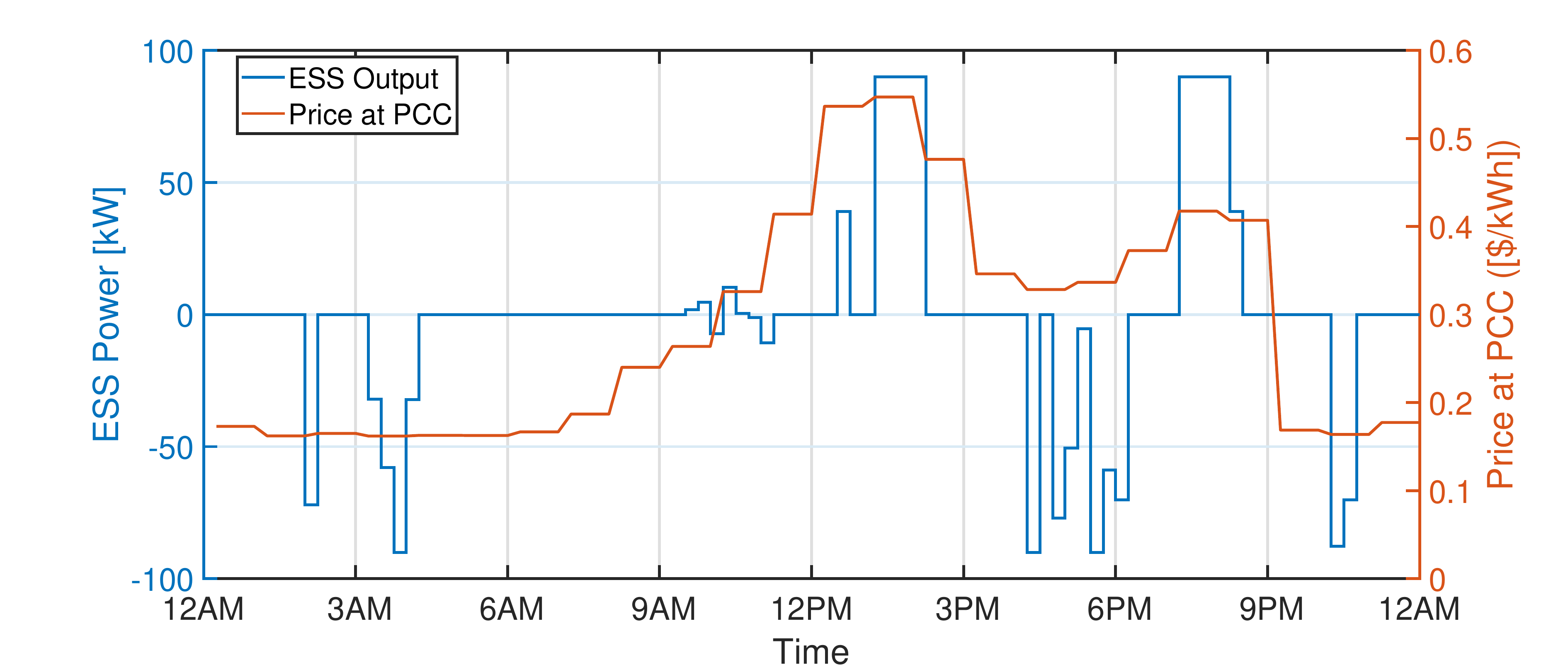}}
	\label{ess_simple}
	\subfloat[Surrogate control]{%
		\includegraphics[width=0.45\textwidth]{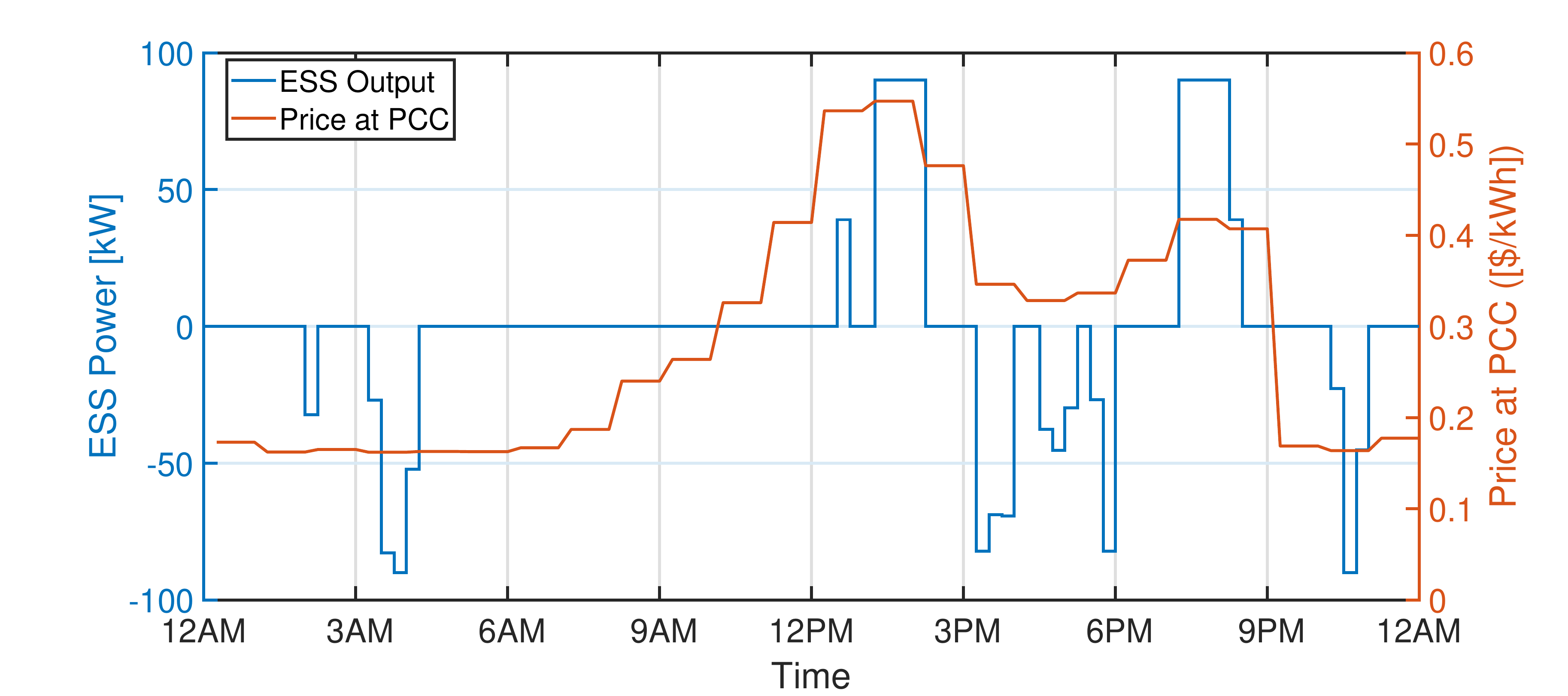}}
	\label{ess_surrogate}
	\caption{The output power of ESS under different control strategy.}
	\label{ess_output}
\end{figure}

\begin{figure}[!htb]
	\centering
	\includegraphics[width=0.45\textwidth]{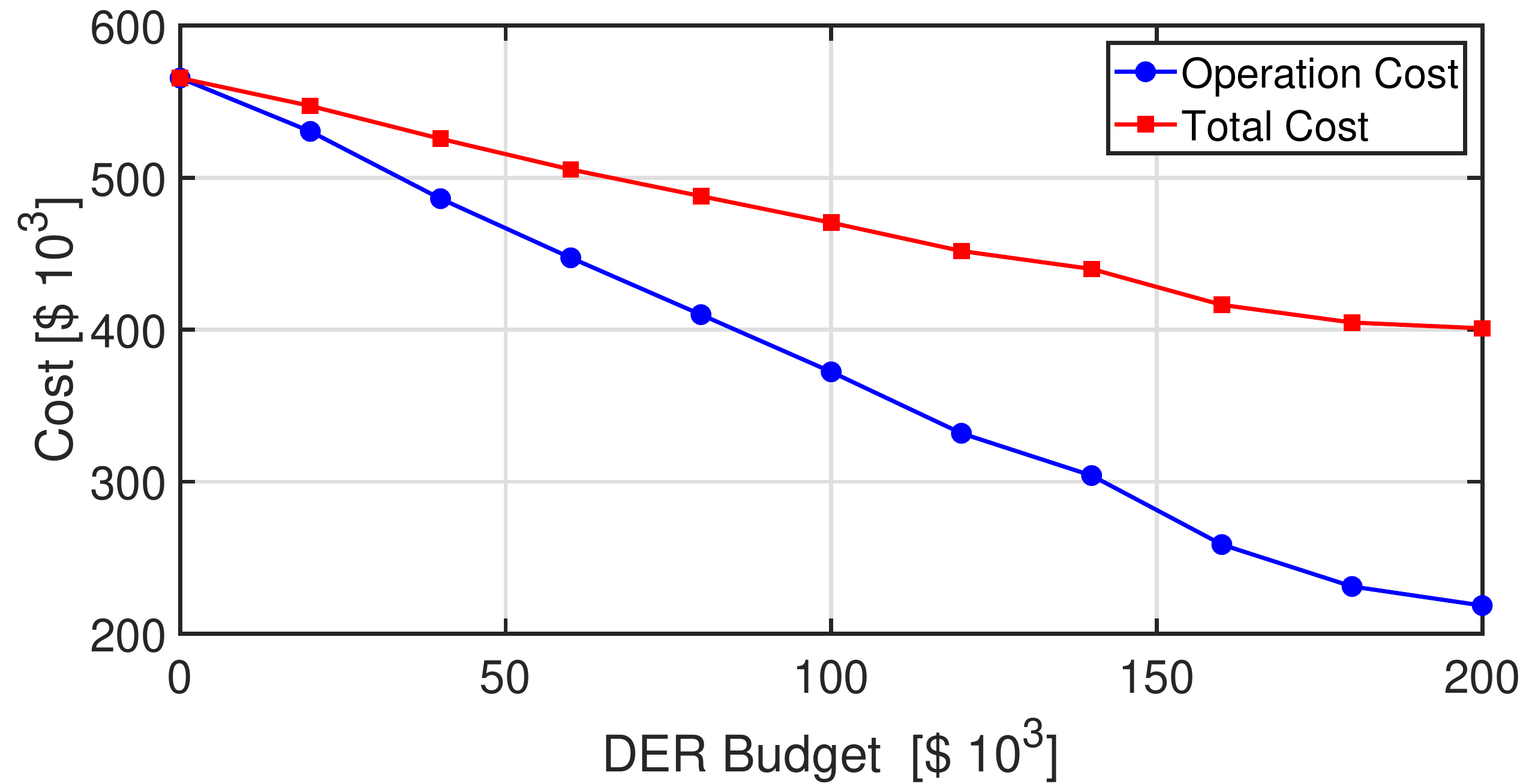}
	\caption{The total and operation cost vs. DER budget.}
	\label{cost_dec}
\end{figure}

Fig. \ref{ess_output} provides the output power of the installed ESS on the Summer day. As can be observed from the figure, the behaviors of ESS under the two control strategies are similar. The ESS will be charged and stored energy during the off-peak hours, and discharge the stored power during the peak hours when the electricity price is high.

Fig. \ref{cost_dec} illustrates the total and operation cost as a function of the budget on DER. Note that the surrogate control strategy is applied for this plot. As can be seen from the figure, both of the operation cost and total cost decrease as the budget on DER increases. However, most of the cost reductions are achieved when the budget on DER is \$180k. The cost can be further reduced if the budget is increased but the amount of the cost reduction is small.

\section{Conclusions and Future Work}
\label{conclusion_future}
This paper proposes an MILP-based optimization model for community MG planning. Different types of DER including ESS, DFG and RES are considered. In addition, the detailed building thermal dynamic characteristics are integrated into the planning model. The optimization model seeks to identify the investment strategy on DER subject to a series of operational and investment constraints. Numerical results based on a community MG consisting of 20 residential buildings demonstrate the effectiveness of the proposed planning model. Moreover, the benefits of including building thermal dynamics in the MG planning model are illustrated.

The islanding capability is an important feature for MG. Future work will extend the planning model to consider the MG islanding uncertainties.

\bibliographystyle{IEEEtran}
\bibliography{IEEEabrv,mybibb}

\end{document}